\documentclass[11pt]{article}

\usepackage{subfigure}
\usepackage{graphicx}
\usepackage{amsmath, amsthm, amssymb}
\usepackage{array}

\usepackage[vcentering,dvips]{geometry}
\newtheorem{theorem}{Theorem}[section]
\newtheorem{lemma}[theorem]{Lemma}

\newtheorem{corollary}[theorem]{Corollary}

\newtheorem{remark}[theorem]{Remark}

\begin{document}

\title{Compressed Sensing with Cross Validation}

\author{Rachel Ward}

\maketitle

\begin{abstract}
Compressed Sensing decoding algorithms can efficiently recover an $N$ dimensional real-valued vector $x$ to within a factor of its best $k$-term approximation by taking $m = 2k\log{N/k}$ measurements $y = \Phi x$.  If the sparsity or approximate sparsity level of $x$ were known, then this theoretical guarantee would imply quality assurance of the resulting compressed sensing estimate.   However, because the underlying sparsity of the signal $x$ is unknown, the quality of a compressed sensing estimate $\hat{x}$ using $m$ measurements is not assured.  Nevertheless, we demonstrate that sharp bounds on the error $ || x - \hat{x} ||_{l^N_2} $ can be achieved with almost no effort.  More precisely, we assume that a maximum number of measurements $m$ is pre-imposed; we reserve $4\log p$ of the original $m$ measurements and compute a sequence of possible estimates  $\big(\hat{x}_j \big)_{j=1}^p$ to $x$ from the $m -4\log{p}$ remaining measurements; the errors $||x - \hat{x}_j ||_{l^N_2} $ for $j = 1, ... , p$ can then be bounded with high probability.  As a consequence, numerical upper and lower bounds on the error between $x$ and the best $k$-term approximation to $x$ can be estimated for $p$ values of $k$ with almost no cost.     Our observation has applications outside of compressed sensing as well.
\\
\\
{\bf Key words:} Compressed Sensing, cross validation, measurements, best $k$-term \\
\indent approximation, Johnson Lindenstrauss Lemma, encoding/decoding, error estimates

\end{abstract}

\section{Introduction}
Compressed Sensing (CS) is a fast developing area in applied mathematics, motivated by the reality that most data we store and transmit contains far less information than its dimension suggests.  For example, a one-dimensional slice through the pixels in a typical grayscale image will contain segments of smoothly varying intensity, with sharp changes between adjacent pixels appearing only at edges in the image. Often this sparsity in information translates into a sparse (or approximately sparse) representation of the data with respect to some standard basis; for the image example, the basis would be a wavelet of curvelet basis.   For such $N$ dimensional data vectors that are well approximated by a $k$-sparse vector (or a vector that contains at most $k << N$ nonzero entries), it is common practice to temporarily store the entire vector, possibly with the intent to go back and replace this vector with a smaller dimensional vector encoding the location and magnitude of its $k$ significant coefficients.  In compressed sensing, one instead collects fewer fixed linear measurements of the data to start with, sufficient in number  to recover the location and numerical value of the $k$ nonzero coordinates at a later time.   Finding "good" linear measurements, as well as fast, accurate, and simple algorithms for recovering the original data from these measurements, are the twofold goals of compressed sensing research today.\\ \\
\noindent { \bf Review of basic CS setup}.  The data of interest is taken to be a real-valued vector $x \in  \mathbb{R}^N$ that is {\it unknown}, but from which we are allowed up to $m < N$ linear measurements, in the form of inner products of $x$ with $m$ vectors $v_j \in \mathbb{R}^N$ of our choosing.  Letting $\Phi$ denote the $m \times N$ matrix whose $j$th row is the vector $v_j$, this is equivalent to saying that we have the freedom to choose and store an $m \times N$ matrix $\Phi$, along with the $m$-dimensional measurement vector $y = \Phi x$.  Of course, since $\Phi$ maps vectors in $\mathbb{R}^N$ to vectors in a smaller dimensional space $\mathbb{R}^m$, the matrix $\Phi$ is not invertible, and we thus have no hope of being able to reconstruct an arbitrary $N$ dimensional vector $x$ from such measurements.\\ \\
\noindent However, if the otherwise unknown vector $x$ is specified to be $k$-sparse, and $k$ is fairly small compared with $N$, then there do exist matrices $\Phi$ for which $y = \Phi x$ uniquely determines $x$, and allows recovery of $x$ using fast and simple algorithms.  It was the interpretation of this phenomenon given by Candes and Tao $\cite{1}$, $\cite{2}$, and Donoho $\cite{3}$, that gave rise to compressed sensing.  In particular, these authors define classes of matrices that possess this property.   One particularly elegant characterization of this class is via the {\it Restricted Isometry Property} (RIP) $\cite{2}$.  A matrix $\Phi$ with unit normed columns is said to be $k$-RIP if all singular values of any $k$ column submatrix of $\Phi$ lie in the interval $[1 - \delta_k, 1 + \delta_k]$ for a given constant $\delta_k < 1$.  
With high probability, $2k$-RIP is obtained, with
\begin{equation}
k = K(m,N) := 2m/ \log(N/m)), \textrm{where } m \leq \frac{1}{2}m,
\label{k}
\end{equation}
for $m \times N$ matrices $\Phi$ whose entries $\Phi_{i,j}$ are independent realizations of a Gaussian or Bernoulli random variable $\cite{4}$.  Also with high probability, an $m \times N$ matrix obtained by selecting $m$ rows at random from the $N \times N$ discrete Fourier matrix satisfies $2k$-RIP of the same order as $\eqref{k}$ up to an additional $\log^3 N$ factor \cite{28}.  In fact, the order of $k$ given by $\eqref{k}$ is optimal given $m$ and $N$, as shown in $\cite{5}$ using classical results on Gelfand widths of $l_1^N$ unit balls in $l_2^N$.   To date, there exist no deterministic constructions of RIP matrices of this order.\\ \\
\noindent { \bf Recovering or approximating x. } As shown in $ \cite{21}$, the following approximation results hold for matrices $\Phi$ that satisfy $2k$-RIP with constant $\delta_{2k} \leq \sqrt{2} - 1$:
\begin{enumerate}
\item If $x \in \mathbb{R}^N$ is $k$-sparse,  then $x$ can be reconstructed from $\Phi$ and the measurement vector $y = \Phi x$ as the solution to the following $\ell_1$ minimization problem:
\begin{equation}
x = {\cal L}_1(\Phi, y) := \arg \min_{\Phi z = y} || z ||_1.
\label{l1}
\end{equation} 
\item If $x$ is not $k$-sparse, the error between $x$ and the approximation $\hat{x} =  {\cal L}_1(\Phi, y)$ is still bounded by
\begin{equation}
|| x - \hat{x} ||_2 \leq \frac{c}{\sqrt{k}} \sigma_k(x)_{l_1^N},
\label{approxerr1}
\end{equation}
where $c = 2(1 + \delta_{2k})/(1 - \delta_{2k})$, and $\sigma_k(x)_{l_p^N} := \inf_{|z| \leq k} ||x - z||_{l_p^N}$ denotes the best possible approximation error in the metric of $l^N_p$ between $x$ and the set of $k$-sparse signals in $\mathbb{R}^N$.  The approximation error $\sigma_k(x)_{l_p^N}$ is realized by the $k$-sparse vector $x_k \in \mathbb{R}^N$ that corresponds to vector $x$ with all but the $k$ largest entries set to zero, independent of the $l^N_p$ norm in the approximation $\sigma_k(x)_{l_p^N}$. 
\end{enumerate} 

\noindent This immediately suggests to use the $\ell_1$-minimizer ${ \cal L}_{1}$ as a means to recover or approximate an unknown $x$ with sparsity constraint.  Several other decoding algorithms are used as alternatives to $\ell_1$ minimization for recovering a sparse vector $x$ from its image $y = \Phi x$, not because they offer better accuracy ( $\ell_1$ minimization gives optimal approximation bounds when $\Phi$ satisfies RIP), but because they can be faster and easier to implement.  For a comprehensive survey on compressed sensing decoding algorithms, we refer the reader to \cite{30}. \\ \\
\noindent{ \bf Estimating the accuracy of CS estimates.}  According to the bound $\eqref{approxerr1}$, the quality of a compressed sensing estimate $\hat{x} = \triangle(\Phi, y)$ depends on how well $x$ can be approximated by a $k$-sparse vector, where the value of $k$ is determined by the number of rows $m$ composing $\Phi$.   While $k$ is assumed to be known and fixed in the compressed sensing literature, no such bound is guaranteed for real-world signal models such as vectors $x \in \mathbb{R}^N$ corresponding to wavelet coefficient sequences of discrete photograph-like images.  Thus,  the quality of a compressed sensing estimate $\hat{x}$ in general is not guaranteed. \\ \\
If the error $||x - \hat{x} ||_2/||x||_2$ incurred by a particular approximation $\hat{x}$ were \emph{observed} to be large, then decoding could be repeated using a larger number of measurements,  perhaps at increasing measurement levels $\{m_1, m_2, ... , m_p \}$, until the error $||x - \hat{x}_j ||_2/||x||_2$ corresponding to $m_j$ measurements were observed to be sufficiently small.  Of course, the errors $||x - \hat{x}_j ||_2$ and $|| x - \hat{x}_j ||_2 / || x ||_2$ are typically not known, as $x$ is unknown.  Our main observation is that one can apply the Johnson-Lindenstrauss lemma $\cite{13}$ to the set of $p$ points,
\begin{equation}
\{ (x - \hat{x}_1), (x - \hat{x}_2), ... , (x - \hat{x}_p) \}.
\end{equation}
In particular, $r = O(\log p)$ measurements of $x$, provided by $y_{\Psi} = \Psi x$, when $\Psi$ is, e.g. a Gaussian or Bernoulli random matrix, are sufficient to guarantee that with high probability,
\begin{equation}
4/5 || y_{\Psi} - \Psi \hat{x}_j ||_2 \leq  ||x - \hat{x}_j ||_2  \leq 4/3 || y_{\Psi} - \Psi \hat{x}_j ||_2 \label{same}
\end{equation}
and
\begin{equation}
1/3 \frac{|| y_{\Psi} - \Psi \hat{x}_j ||_2}{ ||y_{\Psi}||_2} \leq \frac{ ||x - \hat{x}_j ||_2}{ ||x||_2}  \leq 2 \frac{|| y_{\Psi} - \Psi \hat{x}_j ||_2}{||y_{\Psi}||_2}
\label{same2}
\end{equation}
for any $p$ compressed sensing estimates.   The equivalences $\eqref{same}$ and  
$\eqref{same2}$ allow the {\it measurable} quantities $||y_{\Psi} - \Psi \hat{x}_j ||_2$ and $|| y_{\Psi} - \Psi \hat{x}_j ||_2/ ||y_{\Psi}||_2$ to function as proxies for the {\it unknown} quantities $|| x - \hat{x}_j ||_2$ and $||x - \hat{x}_j ||_2 / ||x||_2$; these proxies can be used to
\begin{enumerate}
\renewcommand{\labelenumi}{(\alph{enumi})}
\item provide tight numerical upper and lower bounds on the errors $ || x - \hat{x}_j ||_2$ and $||x - \hat{x}_j ||_2/||x||_2$ at up to $p$ compressed sensing estimates $\hat{x}_j$, 

\item provide estimates of the underlying $k$-term approximations $|| x - x_k ||_2$ of $x$ for up to $p$ different values of $k$, and

\item return from among a sequence of estimates $( \hat{x}_1, ... , \hat{x}_p )$ with different initialization parameters, an estimate $\hat{x}_{cv} = \arg \min_j || y_{\Psi}- \Psi \hat{x}_j ||_2$ having error $|| x - \hat{x}_{cv}||_2$ that does not exceed a small multiplicative factor of the best possible error in the metric of $\ell_2^N$ between $x$ and an element from the sequence at hand.  
\end{enumerate}
More precisely, all CS decoding algorithms require as input a parameter $m$ corresponding to the number of rows in $\Phi$; some compressed sensing decoding algorithms (such as greedy algorithms) require also a parameter $k$ indicating the sparsity level of $x$, and other algorithms require as input a bound $\gamma$ on the expected amount of energy in $x$ outside of its significant coefficients.  All CS decoding algorithms can be symbolically represented by functions of the form $\triangle(\Phi, y, k, \gamma)$,  and we will give examples where each of the parameters $m$, $k$, and $\gamma$ can be optimized over a sequence of estimates $(\hat{x}_1, \hat{x}_2, ..., \hat{x}_p)$ parametrized by increasing hypotheses on each of the variables $m$, $k$, and $\gamma$.  \\ \\
The estimation procedure described above, although novel in its proposed application, is by no means new. {\it Cross validation} is a technique used in statistics and learning theory whereby a data set is separated into a training/estimation set and a test/cross validation set, and the test set is used to prevent overfitting on the training set by estimating underlying noise parameters.  We will take a set of $m$ measurements of $x$, and use $m - r$ of these measurements, $\Phi x$, in a compressed sensing decoding algorithm to return a sequence $(\hat{x}_1, \hat{x}_2, ... )$ of candidate approximations to $x$.  The remaining $r$ measurements, $\Psi x$, are then used to identify from among this sequence a single approximation $\hat{x} = \hat{x}_j$, along with an estimate of the sparsity level of $x$.   The application of cross validation to compressed sensing has been studied by Boufounos, Duarte, and Baraniuk in $\cite{7}$, but in a different context from the present paper, and without the mathematical justification of the Johnson Lindenstrauss lemma that we present below.

 \section{Preliminary Notation}
 Throughout the paper, we will be dealing with large dimensional vectors that have few nonzero coefficients.  We use the notation 
 \begin{equation}
 | x | = n
 \end{equation}
 to indicate that a vector $x \in \mathbb{R}$$^N$ has exactly $n$ nonzero coefficients.\\ \\

\noindent We will sometimes use the notation $a \sim_{\epsilon} b$ as shorthand for the multiplicative relation
\begin{equation} 
(1 - \epsilon)a \leq b \leq (1 + \epsilon)a,
\end{equation}
that can be worded as ``the quantity $a$ approximates the quantity $b$ to within a multiplicative factor of $(1 \pm \epsilon)$".  Note that the relation $\sim_{\epsilon}$ is not symmetric.  Properties of the relation $a \sim_{\epsilon} b$ are listed below; we leave the proofs  (which amount to a string of simple inequalities) as an exercise for the reader.

\begin{lemma} Fix $\epsilon \in (0,1)$.
\begin{enumerate}
\item  If $a, b \in \mathbb{R}^+$ satisfy $a \sim_{\epsilon} b$, then $b/\big[(1+\epsilon)(1-\epsilon)\big] \sim_{\epsilon} a$. 
\item If $a,b,c,d \in \mathbb{R}^+$ satisfy  $a \sim_{\epsilon} b$ and  $c \sim_{\epsilon} d$, then $a/c \sim_{\delta} b/d$ for parameter $\delta = 2\epsilon/1-\epsilon$.
\item If $( a_1, a_2, ... , a_p )$ and $( b_1, b_2, ... , b_p )$ are sequences in $\mathbb{R}^{+}$, and $a_j \sim_{\epsilon} b_j$ for each $1 \leq j \leq p$,  then $\min_j a_j \sim_{\epsilon} \min_j b_j$.
\end{enumerate}
\label{babylemma}
\end{lemma}

\section{Mathematical Foundations}
 
The Johnson Lindenstrauss (JL) lemma, in its original form, states that any set of $p$ points in high dimensional Euclidean space can be embedded into $\epsilon^{-2} \log(p)$ dimensions, without distorting the distance between any two points by more than a factor of $(1 \pm \epsilon)$ $\cite{13}$.   In the same paper, it was shown that a random orthogonal projection would provide such an embedding with positive probability.  Following several simplifications to the original proof $\cite{15}$, $\cite{12}$, $\cite{14}$, it is now understood that Gaussian random matrices, among other purely random matrix constructions, can substitute for the random projection in the original proof of Johnson and Lindenstrauss.  Of the several versions of the lemma now appearing in the literature, the following variant presented in Matousek $\cite{16}$ is most applicable to the current presentation.

\begin{lemma}[Johnson-Lindenstrauss Lemma]
Fix an accuracy parameter $\epsilon \in (0, 1/2]$, a confidence parameter $\delta \in (0,1)$, and an integer $r \geq r_0 = C\epsilon^{-2}\log{\frac{1}{2\delta}}$. \\ \\
Let ${ \cal M}$ be a random $r \times N$ matrix whose entries ${\cal M}_{i,j}$ are independent realizations of a random variable R that satisfies:
\begin{enumerate}
\item $Var(R) = 1/r$  (so that the columns of ${\cal M}$ have expected $\ell_2$ norm 1)
\item $E(R) = 0$,  
\item For some fixed $a > 0$ and for all $\lambda$,
\begin{equation}
\textrm{Prob}[ |R| > \lambda ] \leq 2e^{-a \lambda^2}
\label{a} 
\end{equation}
\end{enumerate}
Then for a predetermined $x \in \mathbb{R}^N$, 
\begin{equation}
(1 - \epsilon)|| x||_{l_2^N} \leq || {\cal M} x ||_{l_2^r} \leq (1 + \epsilon)|| x||_{l_2^N}
\label{jlbound} 
\end{equation}
is satisfied with probability exceeding $1 - \delta$.
\label{psi}
\end{lemma}
\noindent The constant $C$ bounding $r_0$ in Lemma $\eqref{psi}$ grows with the parameter $a$ specific to the construction of ${\cal M}$ $\eqref{a}$.   Gaussian and Bernoulli random variables $R$ will satisfy the concentration inequality $\eqref{a}$ for a relatively small parameter $a$ (as can be verified directly), and for these matrices one can take $C = 8$ in Lemma $\eqref{psi}$. \\ \\
\noindent The Johnson Lindenstrauss lemma can be made intuitive with a few observations.  Since ${\bf E}\big[R\big] =0$ and ${\bf Var}\big[R\big] = \frac{1}{r}$, the random variable $|| {\cal M} x ||_2^2$ equals $||x||_2^2$ in expected value; that is,
\begin{equation}
{\bf E} \big[ \textrm{ } || {\cal M} x ||^2_2\textrm{ }  \big] = || x ||_2^2.
\end{equation} 
Additionally, $|| {\cal M} x ||^2_2$ inherits from the random variable $R$ a nice concentration inequality:
\begin{eqnarray}
\textrm{Prob}\big[ || {\cal M} x ||_2^2 - ||x||_2^2 > \epsilon ||x||_2^2 \big] &\leq& e^{-a(2\epsilon\sqrt{r})^2} 
\leq \delta/2.
\label{bound}
\end{eqnarray} 
The first inequality above is at the heart of the JL lemma; its proof can be found in $\cite{16}$.  The second inequality follows using that $r \geq  (2a\epsilon^2)^{-1}\log(\frac{\delta}{2})$ and $\epsilon \leq 1/2$ by construction.  A bound similar to $\eqref{bound}$ holds for $\textrm{Prob} \big[ || {\cal M} x ||_2^2 - ||x||_2^2 < - \epsilon ||x||_2^2 \big]$ as well, and combining these two bounds gives desired result $\eqref{jlbound}$.\\ \\
\noindent For fixed $x \in \mathbb{R}^N$, a random matrix ${\cal M}$ constructed according to Lemma $\eqref{psi}$ fails to satisfy the concentration bound $\eqref{jlbound}$ with probability at most $\delta$.  Applying Boole's inequality, ${\cal M}$ then fails to satisfy the stated concentration on any of $p$ predetermined points $\{x_j\}_{j=1}^p$, $x_j \in \mathbb{R}^N$, with probability at most $\xi = p\delta$.   In fact, a specific value of $\xi \in (0,1)$ may be imposed for fixed $p$ by setting $\delta = \xi/p$.  These observations are summarized in the following corollary to Lemma $\eqref{psi}$.

\begin{corollary}
Fix an accuracy parameter $\epsilon \in (0,1/2]$, a confidence parameter $\xi \in (0,1)$, and fix a set of $p$ points $\{x_j\}_{j=1}^p \subset \mathbb{R}^N$.  Set $\delta = \xi/p$, and fix an integer $r \geq r_0 = C\epsilon^{-2}\log{\frac{1}{2\delta}} = C\epsilon^{-2}\log{\frac{p}{2\xi}}$.  If ${\cal M}$ is a $r \times N$ matrix constructed according to Lemma $\eqref{psi}$, then with probability $\geq 1 - \xi$, the bound
\begin{equation}
(1 - \epsilon) || x_j ||_{l_2^N} \leq || {\cal M} x_j ||_{l_2^r} \leq (1 + \epsilon) || x_j ||_{l_2^N}
\label{jlplus}
 \end{equation}
obtains for each $j=1,2, ... ,p$.
\label{myjl}
\end{corollary}

\section{Cross Validation in Compressed Sensing}
We return to the situation where we would like to approximate a vector $x \in \mathbb{R}^N$ with an assumed sparsity constraint using $m < N$ linear measurements $y = {\cal A} x$ where ${\cal A}$ is an $m \times N$ matrix of our choosing.  Continuing the discussion in Section 1, we will not reconstruct $x$ in the standard way by $\hat{x} = \Delta({\cal A}, y, k, \gamma)$ for fixed values of the input parameters, but instead separate the $m \times N$ matrix ${\cal A}$ into an $n \times N$ {\it implementation} matrix $\Phi$ and an $r \times N$ {\it cross validation} matrix $\Psi$, and separate the measurements $y$ accordingly into $y_{\Phi}$ and $y_{\Psi}$.  We use the implementation matrix $\Phi$ and corresponding measurements $y_{\Phi}$ as input into the decoding algorithm to obtain a sequence of possible estimates $(\hat{x}_1, ..., \hat{x}_p)$ corresponding to increasing one of the input parameters $m$, $k$, or $\gamma$.  We reserve the cross validation matrix $ \Psi$ and measurements $y_{\Psi}$ to estimate each of the error terms $|| x - \hat{x}_j ||_2$ in terms of the computable $||y_{\Psi} - \Psi \hat{x}_j||_2$. Our main result, which follows from Corollary $\eqref{myjl}$, details how the number of cross validation measurements $r$ should be chosen in terms of the desired accuracy $\epsilon$ of estimation, confidence level $\xi$ in the prediction, and number $p$ of estimates $\hat{x}_j$ to be measured: 
\begin{theorem}
For a given accuracy $\epsilon \in (0,1/2]$, confidence $\xi \in (0,1)$, and number $p$ of estimates $\hat{x}_j \in \mathbb{R}^N$, it suffices to allocate $r  = \lceil C\epsilon^{-2}\log{\frac{p}{2\xi}} \rceil$ rows to a cross validation matrix $\Psi$ of Gaussian or Bernoulli type, normalized according to Lemma \eqref{myjl} and independent of the estimates $\hat{x}_j$, to obtain with probability greater than or equal to $1 - \xi$, and for each $j = 1,2, ..., p$, the bounds
\\
\\
\begin{equation}
 \frac{1}{1+\epsilon} \leq \frac{||x - \hat{x}_j ||_{l_2^N}}{|| \Psi (x - \hat{x}_j) ||_{l_2^{\ell}}}\leq \frac{1}{1-\epsilon}
 \label{1a}
\end{equation}
and 
\begin{eqnarray}
 \frac{1-3\epsilon}{(1+\epsilon)(1-\epsilon)^2} &\leq& \frac{||x - \hat{x}_j ||_{l_2^N}/||x||_2}{|| \Psi (x - \hat{x}_j) ||_{l_2^{\ell}}/||\Psi x ||_2} \leq \frac{1}{(1-\epsilon)^2} \nonumber \\ \nonumber \\
\label{2a}
\end{eqnarray}
and also
\begin{equation}
 \frac{1}{1+\epsilon} \leq \frac{\eta_{or}}{\widehat{\eta_{cv}}} \leq \frac{1}{1-\epsilon}
\label{1b}
\end{equation}
where $\eta_{or} = \min_{1\leq j \leq p}  ||x - \hat{x}_j ||_{l_2^N}$ is the unknown  \emph{oracle error} corresponding to the best possible approximation to $x$ in the metric of $l_2^N$ from the sequence $(\hat{x}_1, ..., \hat{x}_p)$, and $\widehat{\eta_{cv}} = \min_{1 \leq j \leq p} || \Psi (x - \hat{x}_j) ||_{l_2^r}$ is the observable cross validation error.  
\label{mainthm}
\end{theorem} 

\begin{proof}.
\begin{itemize}
\item The bounds in \eqref{1a} are obtained by application of Lemma $\eqref{myjl}$ to the $p$ points $u_j = x - \hat{x}_j$, and rearranging the resulting bounds according to Lemma $\eqref{babylemma}$ part $(1)$.   The bound \eqref{1b} follows from the bounds \eqref{1a} and part (3) of Lemma $\eqref{babylemma}$.
\item  The bounds in \eqref{2a} are obtained by application of Lemma $\eqref{myjl}$ to the $p+1$ points $u_0 = x, u_j = x - \hat{x}_j$, and regrouping the resulting bounds according to part (2) of Lemma $\eqref{babylemma}$.
\end{itemize}
\end{proof}
\begin{remark}
\emph{The measurements making up the cross validation matrix $\Psi$ must be independent of the measurements comprising the rows of the implementation matrix $\Phi$.  This comes from the requirement in Lemma $\eqref{psi}$ that the matrix $\Psi$ be independent of the points $u_j = x - \hat{x}_j$.  This requirement is crucial, as observed  when $\hat{x}$ solves the $\ell_1$ minimization problem 
\begin{equation}
\hat{x} = \arg \min_{z \in \mathbb{R}^N} || z ||_1 \textrm{ subject to } \Phi z = \Phi x,
\end{equation}
in which case the constraint $\Phi(\hat{x} - x) = 0$ clearly precludes the rows of $\Phi$ from giving any information about the error $|| \hat{x} - x ||_2$. }
\end{remark}

\begin{remark}
\emph{Theorem $\eqref{mainthm}$ should be applied with a different level of care depending on what information about the sequence $(x - x_1, x - x_2, ..., x - x_p)$ is sought.  If the minimizer $\hat{x} = \arg \min_{1 \leq j \leq p} || \Psi (x - \hat{x}_j) ||_{l_2^r}$ is sufficient for one's purposes, then the precise normalization of $\Psi$ in Theorem $\eqref{mainthm}$ is not important.  The normalization doesn't matter either for estimating the normalized quantities $||x - \hat{x}_j||_2/||x||_2$.  On the other hand, if one is using cross validation to obtain estimates for the quantities $||x - \hat{x}_j||_2$, then normalization is absolutely crucial, and one must observe the normalization factor given by Lemma $\eqref{myjl}$ that depends on the number of rows $r$ allocated to the cross validation matrix $\Psi$.}
\end{remark}

\section{Applications of cross validation to compressed sensing}

\subsection{Estimation of the best $k$-term approximation error}
We have already seen that if the $m \times N$ matrix $\Phi$ satisfies $2k$-RIP with parameter $\delta \leq \sqrt{2}-1$, and $\hat{x} = {\cal L}_1(\Phi, \Phi x)$ is returned as the solution to the $\ell_1$ minimization problem $\eqref{l1}$, then the error between $x$ and the approximation $\hat{x}$ is bounded by
\begin{equation}
|| x - \hat{x} ||_2 \leq \frac{c}{\sqrt{k}} \sigma_k(x)_{l_1^N}.
\label{approxerr}
\end{equation}
Several other decoding algorithms in addition to $\ell_1$ minimization enjoy the reconstruction guarantee $\eqref{approxerr}$ under similar bounds on $\Phi$, such as the Iteratively Reweighted Least Squares algorithm (IRLS) $\cite{29}$, and the greedy algorithms CoSAMP $\cite{30}$ and Subspace Pursuit $\cite{31}$.  It has recently been shown \cite{18} \cite{20} that if the bound \eqref{approxerr} is obtained, and if $x - \hat{x}$ lies in the null space of $\Phi$ (as is the case for the decoding algorithms just mentioned), then if $\Phi$ is a Gaussian or a Bernoulli random matrix, the error $||x - \hat{x} ||_2$ also satisfies a bound, with high probability on $\Phi$, with respect to the $\ell_2^N$ residual, namely
\begin{equation}
|| x - \hat{x} ||_2 \leq c\sigma_k(x)_{l_2^N},
\label{w2}
\end{equation}
for a reasonable constant $c$ depending on the RIP constant $\delta_{2k}$ of $\Phi$.  In the event that $\eqref{w2}$ is obtained, a cross validation estimate $|| \Psi (x - \hat{x}) ||_{l_2^{r}}$ can be used to lower bound the residual $\sigma_k(x)_{l_2^N}$, with high probability, according to
\begin{eqnarray}
(1 - \epsilon) || \Psi (x - \hat{x}) ||_{l_2^{\ell}} \leq ||x - \hat{x} ||_{l_2^N}  \leq c \sigma_k(x)_{l_2^N},
\end{eqnarray}
with $O(\frac{1}{\epsilon^2})$ rows reserved for the matrix $\Psi$ $\eqref{mainthm}$.  At this point, we will use Corollary 3.2 of $\cite{8}$, where it is proved that if the bound $\eqref{approxerr}$ holds for $\hat{x}$ with constant $c$, then the same bound will hold for
\begin{equation}
\hat{x}_k = \arg \min_{z: |z| \leq k} || \hat{x} - z ||_{l_2^N},
\end{equation}
the best $k$-sparse approximation to $\hat{x}$, with constant $\tilde{c} = 3c$.  Thus, we may assume without loss of generality that $\hat{x}$ is $k$-sparse,  in which case $|| \Psi (x - \hat{x}) ||_{l_2^{r}}$ also provides an upper bound on the residual $\sigma_k(x)_{l_2^N}$ by
\begin{equation}
(1 + \epsilon) || \Psi (x - \hat{x}) ||_{l_2^{r}}  \geq  ||x - \hat{x} ||_{l_2^N} \geq \sigma_k(x)_{l_2^N}.
\end{equation} 
With almost no effort then, cross validation can be incorporated into many decoding algorithms to obtain tight upper and lower bounds on the unknown $k$-sparse approximation error $\sigma_k(x)_{l_2^N}$ of $x$.  More generally, the allocation of $10\log{p}$ measurements to the cross validation matrix $\Psi$ is sufficient to estimate the errors $(||x - x_{k_j}||_2)_{j=1}^p$ or the normalized approximation errors $(||x - x_{k_j}||_2/||x||_2)^p_{j=1}$ at $p$ sparsity levels $k_j$ by decoding $p$ times, adding $m_j$ measurements to the implementation matrix $\Phi_j$ at each repetition.  Recall that the quantities $k_j$ and $m_j$ are related by $k_j = 2m_j/ \log(N/m_j))$ according to $\eqref{k}$. 

\subsection{Choice of the number of measurements $m$}
Photograph-like images have wavelet or curvelet coefficient sequences $x \in \mathbb{R}^N$ that are \emph{compressible} $\cite{23}$ $\cite{32}$, having entries that obey a power law decay
\begin{equation}
|x|_{(k)} \leq c_s k^{-s},
\label{power}
\end{equation}
where $x_{(k)}$ denotes the $k$th largest coefficient of $x$ in absolute value, the parameter $s > 1$ indicates the level of compressibility of the underlying image, and $c_s$ is a constant that depends only on $s$ and the normalization of $x$.  From the definition $\eqref{power}$, compressible signals are immediately seen to satisfy
\begin{equation}
||x - x_k ||_1/\sqrt{k} \leq c_s' k^{-s + 1/2}, 
\label{compress}
\end{equation}
so that the solution $\hat{x}_m = {\cal L}_1(\Phi, \Phi x)$ to the $\ell_1$ minimization problem $\eqref{l1}$ using an $m \times N$ matrix $\Phi$ of optimal RIP order $k = 2m/ \log(N/m))$ satisfies 
\begin{equation}
|| x - \hat{x}_{m} ||_2 \leq c_{s,\delta} k^{-s+1/2}.
\label{levelcomp}
\end{equation}
The number of measurements $m$ needed to obtain an estimate $\hat{x}_{m}$ satisfying $ ||x - \hat{x}_m ||_2 \leq \tau$ for a predetermined threshold $\tau$ will vary according to the compressibility of the image at hand.  Armed with a total of $m$ measurements, the following decoding method that \emph{adaptively} chooses the number of measurements for a given signal $x$ presents a more democratic alternative to standard compressed sensing decoding structure:  

\begin{table}[h]
\label{BDC}
{ \small
\begin{center}
\caption{ {\it \small CS decoding structure with adaptive number of measurements} } 
\end{center}
\begin{enumerate}
\item { \it Input}: The $m$-dimensional vector $y = \Phi x$, the $m \times N$ matrix $\Phi$, (in some algorithms) the sparsity level $k$, and  (again, in some algorithms) a bound $\gamma$ on the noise level of $x$, the number $p$ of of row subsets of $\Phi$,  $(\Phi_1, \Phi_2, ..., \Phi_p)$, corresponding to increasing number of rows $m_1 < m_2 < ... < m_p < m$, and threshold $\tau > 0$.
\item { \it Initialize} the decoding algorithm at $j = 1$.
\item {\it Estimate} $\hat{x}_j = \triangle(\Phi_{m_j}, y_{m_j}, k, \gamma)$ with the decoder $\triangle$ at hand, using only the first $m_j$ measurement rows of $\Phi$.  The previous estimate $\hat{x}_{j-1}$ can be used for ``warm initialization" of the algorithm, if applicable.  The remaining $r_j = m - m_j$ measurement rows are allocated to a cross validation matrix $\Psi_j$ that is used to estimate the resulting error $|| x - \hat{x}_j ||_2/||x||_2$. 
\item { \it Increment}  $j$ by 1, and iterate from step 3 if stopping rule 
is not satisfied.
\item { \it Stop}: at index $j = j^* < p $ if $||x - x_{m_j}||_2/|| x ||_2 \leq \tau$ holds with near certainty, as indicated by
\begin{equation}
\frac{\sqrt{r_j}||\Psi (x - x_{m_j}) ||_2/ || \Phi x||_2}{\sqrt{r_j} - 3\log{p}} \leq \tau
\label{r_j}
\end{equation}
according to Theorem $\eqref{mainthm}$.  If the maximal number of decoding measurements $m_p < m$ have been used at iteration $p$, and $\eqref{r_j}$ indicates that $|| x - \hat{x}_{m_p} ||_2/||x||_2 > \tau$ still, return $\hat{x}_m = \triangle(\Phi, y, k, \gamma)$ using all $m$ measurements, but with a warning that the underlying image $x$ is probably too dense, and its reconstruction is not trustworthy.     
\label{CS decode}
\end{enumerate}  
}
\end{table}

\subsection{Choice of regularization parameter in homotopy-type  algorithms}
Certain compressed sensing decoding algorithms iterate through a sequence of intermediate estimates $\hat{x}_j$ that could be potential optimal solutions to $x$ under certain reconstruction parameter choices.  This is the case for greedy and homotopy-continuation based algorithms.  In this section, we study the application of cross validation to the intermediate estimates of decoding algorithms of homotopy-continuation type.  
\\
\\
LASSO is the name coined in $\eqref{33}$ for the problem of minimizing of the following convex program:
\begin{equation}
\hat{x}^{[\tau]} = \arg \min_{z \in \mathbb{R}^N} || \Phi x - \Phi z ||_{\ell_{2}} + \tau ||z||_1
\label{lasso}
\end{equation}
The two terms in the LASSO optimization problem $\eqref{lasso}$ enforce data fidelity and sparsity, respectively, as balanced by the regularization parameter $\tau$.  In general, choosing an appropriate value for $\tau$ in $\eqref{lasso}$ is a hard problem; when $\Phi$ is an underdetermined matrix, as is the case in compressed sensing, the function $f(\tau) = || x - \hat{x}^{[\tau]} ||_2$ is unknown to the user but is seen empirically to have a minimum at a value of $\tau$ in the interval $[0, || \Phi x ||_{\infty}]$ that depends on the unknown noise level and/or and compressibility level of $x$.      
\\
\\
The \emph{homotopy continuation} algorithm \cite{26}, which can be viewed as the appropriate variant of LARS \cite{26}, is one of many algorithms for solving the LASSO problem $\eqref{lasso}$ at a predetermined value of $\tau$; it proceeds by first initializing $\tau' $ to a value sufficiently large to ensure that the $\ell_1$ penalization term in $\eqref{lasso}$ completely dominates the minimization problem and $x^{[\tau']} = 0$ trivially.  The homotopy continuation algorithm goes on to generate $x^{[\tau']} = 0$ for decreasing $\tau'$ until the desired level for $\tau$ is reached.   If $\tau = 0$, then the homotopy method traces through the entire solution path $x^{[\tau']} \in \mathbb{R}^N$ for $\tau' \geq 0$ before reaching the final algorithm output $x^{[0]} = {\cal L}_1(\Phi, y)$ corresponding to the $\ell_1$ minimizer $\eqref{l1}$.   
\\
\\
From the non-smooth optimality conditions for the convex functional $\eqref{lasso}$, it can be shown that the solution path $\hat{x}^{[\tau]} \in \mathbb{R}^N$ is a piecewise-affine function of $\tau$ \cite{26}, with ``kinks" possible only at a finite number of points $\tau \in \{\tau_1, \tau_2, ... \}$.  Theorem $\eqref{mainthm}$ suggests a method whereby an appropriate value of $\tau*$ can be chosen from among a subsequence of the kinks $(\tau_1, \tau_2, ... , \tau_p)$ by solving the minimization problem $\hat{x}^{[\tau*]} = \arg \min_{j \leq p} || \Psi (x - \hat{x}^{[\tau_j]}) ||_2$ for appropriate cross validation matrix $\Psi$.  Moreover, since the solution $x - \hat{x}^{[\tau]}$ for $\tau_j \leq \tau \leq \tau_{j+1}$ is restricted to lie in the two-dimensional subspace spanned by $x - \hat{x}^{[\tau_j]}$ and $x - \hat{x}^{[\tau_{j+1}]}$,  one can combine the Johnson Lindenstrauss Lemma with a covering argument analogous to that used to derive the RIP property for Gaussian and Bernoulli random matrices in $\cite{4}$, to cross validate the entire \emph{continuum} of solutions $\hat{x}^{[\tau]}$ between $\tau_1 \leq \tau \leq \tau_p$.  More precisely, the following bound holds under the conditions outlined in Theorem $\eqref{mainthm}$, with the exception that $2r$ (as opposed to $r$) measurements are reserved to $\Psi$:

\begin{equation}
 \frac{1}{1+\epsilon} \leq \frac{\min_{\tau_1 \leq \tau \leq \tau_p} || x - \hat{x}^{[\tau]} ||_2}{\min_{\tau_1 \leq \tau \leq \tau_p} ||\Psi ( x - \hat{x}^{[\tau]})||_2} \leq \frac{1}{1-\epsilon}
 \label{lassobound}
\end{equation}

Unfortunately, the bound $\eqref{lassobound}$ is not strong enough to \emph{provably} evaluate the entire solution path $\hat{x}^{[\tau]}$ for $\tau \geq 0$, because the best upper bound on the number of kinks on a generic  LASSO solution path can be very large.  One can prove that this number is bounded by $3^N$, by observing that if $\hat{x}^{[\tau_1]}$ and $\hat{x}^{[\tau_2]}$ have the same sign pattern, then $\hat{x}^{[\tau]}$ also has the same sign pattern for $\tau_1 \leq \tau \leq \tau_2$.  Applying Theorem $\eqref{mainthm}$ to $p = 3^N$ points $x - \hat{x}_j$, this suggests that $O(N)$ rows would need to be allocated to a cross validation matrix $\Psi$ in order for Theorem $\eqref{mainthm}$ and the corollary $\eqref{lassobound}$ to apply to the entire solution path, which clearly defeats the compressed sensing purpose.  However, whenever the matrix $\Phi$ is an $m \times N$ compressed sensing matrix of random Gaussian, Bernoulli, or partial Fourier construction, it is observed \emph{empirically} that the number of kinks along a homotopy solution path is bounded by $3m$, independent of the underlying vector $x \in \mathbb{R}^N$ used to generate the path.  This suggests, at least heuristically, that the allocation of $O(\log{m})$ out of $m$ compressed sensing measurements of this type suffices to ensure that the error $|| x - \hat{x}^{[\tau]} ||_2$ for the solution $\hat{x}^{[\tau]} = \arg \min_{\tau \geq 0} || \Psi (x - \hat{x}^{[\tau]}) ||_2$ will be within a small multiplicative factor of the best possible error in the metric of $\ell_2^N$ obtainable by any approximant $\hat{x}^{[\tau]}$ along the solution curve $\tau \geq 0$.   At the value of $\tau$ corresponding to $\hat{x}^{[\tau]}$, the LASSO solution $\eqref{lasso}$ can be computed using all $m$ measurements $\Phi$ as a final approximation to $x$.
\\
\\
The \emph{Dantzig selector} (DS) \cite{22} refers to a minimization problem that is similar in form to the LASSO problem: 
\begin{equation}
\hat{x}_{\tau} = \arg \min_{z \in \mathbb{R}^N} || \Phi x - \Phi z ||_{\ell_{\infty}} + \tau || z ||_1
\label{DS} 
\end{equation}
The difference between the DS $\eqref{DS}$ and LASSO $\eqref{lasso}$ is the choice of norm ($\ell_{\infty}$ versus $\ell_2$) on the fidelity-promoting term.  Homotopy-continuation based algorithms have also been developed to solve the minimization problem $\eqref{DS}$ by tracing through the solution path $\hat{x}_{\tau'}$ for $\tau' \geq \tau$.  As the minimization problem $\eqref{DS}$ can be reformulated as a linear program, its solution path $\hat{x}_{\tau} \in \mathbb{R}^N$ is seen to be a piecewise \emph{constant} function of $\tau$, in contrast to the LASSO solution path.  In practice, the total number of breakpoints $(\tau_1, \tau_2, ... )$ in the domain $0 \leq \tau$ is observed to be on the same order of magnitude as $m$ when the $m \times N$ matrix $\Phi$ satisfies RIP \cite{24}; thus, the procedure just described to cross validate the LASSO solution path can be adapted to cross validate the solution path of $\eqref{DS}$ as well. 
\\
\\
Thus far we have not discussed the possibility of using cross validation as a stopping criterion for homotopy-type decoding algorithms.  Along the LARS homotopy curve $\eqref{lasso}$, most of the breakpoints $(\tau_1, \tau_2, ... )$ appear only near the end of the curve in a very small neighborhood of $\tau = 0$.  These breakpoints incur only miniscule changes in the error $|| x - \hat{x}_{\tau_j} ||_2$ even though they account for most of the computational expense of the LARS decoding algorithm.  Therefore, it would be interesting to adapt such algorithms, perhaps using cross validation, to stop once $\tau^*$ is reached for which the error $||x - \hat{x}_{\tau^*} ||_2$ is sensed to be sufficiently small.

\subsection{Choice of sparsity parameter in greedy-type algorithms}
Greedy compressed sensing decoding algorithms also iterate through a sequence of intermediate estimates $\hat{x}_j$ that could be potential optimal solutions to $x$ under certain reconstruction parameter choices.  \emph{Orthogonal Matching Pursuit} (OMP), which can be viewed as the prototypical greedy algorithm in compressed sensing, picks columns from the implementation matrix $\Phi$ one at a time in a greedy fashion until, after $k$ iterations, the $k$-sparse vector $\hat{x}_k$, a linear combination of the $k$ columns of $\Phi$ chosen in the successive iteration steps, is returned as an approximation to $x$.  The OMP algorithm is listed in Table 2.  Although we will not describe the algorithm in full detail, a comprehensive study of OMP can be found in $\cite{6}$.   
\begin{table}[h]
{ \small 
\label{omp}
 \begin{center}
\caption{{\it \small Orthogonal Matching Pursuit Basic Structure}}
\end{center}
\begin{enumerate}
\item { \it Input}: The $m$-dimensional vector $y = {\cal B} x$, the $m \times N$ encoding matrix $\Phi$ whose $j^{th}$ column is labeled $\phi_j$, and the sparsity bound $k$.
\item { \it Initialize} the decoding algorithm at $j = 1$, the residual $r_0 = y$, and the index set $\Lambda_0 = \emptyset$.
\item {\it Estimate} 
\begin{enumerate}
\item Find an index $\lambda_j$ that realizes the bound $ (\Phi^T r_{j-1} )_{\lambda_j} = || \Phi^T r_{j-1} ||_{\infty} $.
\item Update the index set $\Lambda_j = \Lambda_{j-1} \cup { \lambda_j}$ and the submatrix of contributing columns: $\Phi_j = [ \Phi_{j-1}\textrm{,  }  \phi_{\lambda_j} ]$
\item Update the residual:
\begin{eqnarray}
s_j &=& \arg \min_x || \Phi_j x - y ||_2 = (\Phi_j^T\Phi_j)^{-1}\Phi_j^Ty, \nonumber \\
a_j &=& \Phi_j x_j \nonumber \\
r_j &=& r_{j-1} - a_j. \nonumber
\end{eqnarray}
\item The estimate $\hat{x}_j$ for the signal has nonzero indices at the components listed in $\Lambda_j$, and the value of the estimate $\hat{x}_j$ in component $\lambda_i$ equals the $ith$ component of $s_j$.
\end{enumerate} 
\item { \it Increment}  $j$ by 1 and iterate from step 3, if $j < k$.
\item { \it Stop}: at $j = k$.  Output $\hat{x}_{omp} = \hat{x}_k$ as approximation to $x$.
\end{enumerate}
}
\end{table}
Note in particular that OMP requires as input a parameter $k$ corresponding to the expected sparsity level for $x \in \mathbb{R}^N$.  Such input is typical among greedy algorithms in compressive sensing (in particular, we refer the reader to $\cite{30}$,  $\cite{29}$, and  $\cite{31}$).  As shown in $\cite{6}$, OMP will recover with high probability a vector $x$ having at most $k \leq m / log(N)$ nonzero coordinates from its image $\Phi x$ if $\Phi$ is a (known) $m \times N$ Gaussian or Bernoulli matrix with high probability.   Over the more general class of vectors $x = x_d + {\cal N}$ that can be decomposed into a $d$-sparse vector $x_d$ (with $d$ presumably less than or equal to $k$) and additive noise vector ${\cal N}$, we might expect an intermediate estimate $\hat{x}_s$ to be a better estimate to $x$ than the final OMP output $\hat{x}_k$, at least when  $d << k$.  Assuming that the signal $x$ admits a decomposition of the form $x = x_d + {\cal N}$, the sequence of intermediate estimates $(\hat{x}_1, ... , \hat{x}_k)$ of an OMP algorithm can be cross validated in order to estimate the noise level and recover a better approximation to $x$.  We will study this particular application of cross validation in more detail below.

\section{Orthogonal Matching Pursuit: A case study}
As detailed in Table 2, a single index $\lambda_j$ is added to a set $\Lambda_j$ estimated as the $j$ most significant coefficients of $x$ at each iteration $j$ of OMP; following the selection of $\Lambda_j$, an estimate $\hat{x}_j$ to $x$ is determined by the least squares solution,
\begin{equation}
\hat{x}_j = \arg \min_{\textrm{supp}(z) \in \Lambda_j} || \Phi z - y ||_2,
\end{equation}
among the subspace of vectors $z \in \mathbb{R}$$^N$ having nonzero coordinates in the index set $\Lambda_j$.  OMP continues as such, adding a single index $\lambda_j$ to the set $\Lambda_j$ at iteration $j$, until $j = k$ at which point the algorithm terminates and returns the $k$-sparse vector $\hat{x}_{omp} = \hat{x}_k$ as approximation to $x$.  
\\
\\
Suppose $x$ has only $d$ significant coordinates.  If $d$ could be specified beforehand, then the estimate $\hat{x}_d$ at iteration $j = d$ of OMP would be returned as an approximation to $x$.   However, the sparsity $d$ is not known in advance, and $k$ will instead be an upper bound on $d$.  As the estimate $\hat{x}_j$ in OMP can be then identified with the hypothesis that $x$ has $j$ significant coordinates, the application of cross-validation as described in the previous section applies in a very natural way to OMP.  In particular, we expect $\hat{x}_{or}$ and $\hat{x}_{cv}$ of Theorem $\eqref{mainthm}$ to be close to the estimate $\hat{x}_j$ at index $j = |x|$ corresponding to the true sparsity of $x$; furthermore, in the case that $|x|$ is significantly less than $k$, we expect the cross validation estimate $\hat{x}_{cv}$ to be a better approximation to $x$ than the OMP-returned estimate $\hat{x}_k$.  We will put this intuition to the test in the following numerical experiment. 

\subsection{Experimental setup}
We initialize a signal $x_0$ of length $N = 3600$ and sparsity level $d = 100$ as
\begin{equation}
 x_0(j) = 
 \left \{ \begin{array}{ll}
1,  &  \textrm{ for  } j = 1... 100 \\
0, & \textrm{ else}.
\end{array} \right.
\label{x0}
\end{equation}
\noindent Noise is then added to $x_a = x_0 + {\cal N}_a$ in the form of  a Gaussian random variable ${\cal N}_a$ distributed according to
\begin{equation}
{\cal N}_a \sim N(0, .05),
\end{equation}
and the resulting vector $x_a$ is renormalized to satisfy $||x_a||_{l_2^N} = 1$.  This yields an expected noise level of 
\begin{equation}
E(\sigma_d(x_a)) \approx .284.
\end{equation}
\\
We fix the input $k=200$ in Table 2, and assume we have a total number of compressed sensing measurements $m = 800$.  A number $r$ of these $m$ measurements are allotted to cross validation, while the remaining $n = m - r$ measurements are allocated as input to the OMP algorithm in Table 2. This experiment aims to numerically verify Theorem $\eqref{mainthm}$; to this end, we specify a confidence $\xi = 1/100$, and solve for the accuracy $\epsilon$ according to the relation $r =  \epsilon^{-2} \log(\frac{k}{2\xi})$; that is, 

\begin{equation}
\epsilon(r) = \sqrt{\frac{\log(\frac{k}{2\xi})}{r}} \approx \frac{3}{\sqrt{r}}.
\label{rach}
\end{equation} 

\noindent Note that the specification $\eqref{rach}$ corresponds to setting the constant $C = 1$ in Theorem $\eqref{mainthm}$.  Although $C \geq 8$ is needed for  the proof of the Johnson Lindenstrauss lemma at present, we find that in practice $C = 1$ already upper bounds the optimal constant needed for Theorem $\eqref{mainthm}$ for Gaussian and Bernoulli random ensembles. 
\\
\\
A single (properly normalized) Gaussian $n \times N$ measurement matrix $\Phi$ is generated (recall that $n$ = $m$ - $r$) , and this matrix and the measurements $y = \Phi x$ are provided as input to the OMP algorithm; the resulting sequence of estimates  $(\hat{x}_1, \hat{x}_2, ..., \hat{x}_k)$ is stored.  The final estimate $\hat{x}_k$ from this sequence is the returned OMP estimate $\hat{x}_{omp}$ to $x$.  The error $\eta_{omp} = || \hat{x}_{omp} - x ||_2$  is greater than or equal to the oracle error of the sequence, $\eta_{or} = \min_{\hat{x}_j} || x - \hat{x}_j ||_2$.
\\
\\
With the sequence $(\hat{x}_1, \hat{x}_2, ... , \hat{x}_k)$ at hand, we consider $1000$ realizations $\Psi_q$ of an $r \times N$ cross validation matrix having the same componentwise distribution as $\Phi$, but normalized to have variance $1/r$ according to Theorem $\eqref{psi}$.  The cross validation error
 \begin{equation}
 \widehat{\eta_{cv}}(q) = \min_j ||\Psi_q (x - \hat{x}_j)||_{l_2^r}
 \end{equation}
 is measured at each realization $\Psi_q$; we plot the average $\bar{\widehat{\eta_{cv}}}$ of these $1000$ values and intervals centered at $\bar{\widehat{\eta_{cv}}}$ having length equal to twice the empirical standard deviation.  Note that we are effectively testing $1000$ trials of OMP-CV, the algorithm which modifies OMP to incorporate cross validation so that $(\hat{x}_{cv}, \widehat{\eta_{cv}})$ are output instead of $\hat{x}_{omp} = \hat{x}_k$.
 \\
 \\
 At the specified value of $\xi$, Theorem $\eqref{mainthm}$ part $\eqref{2a}$ (with constant $C = 1$) implies that 
 \begin{equation}
 \big(1 -\epsilon)\eta_{or} \leq \widehat{\eta_{cv}}(q) \leq \big( 1 +\epsilon \big)\eta_{or}
 \label{theory}
 \end{equation}
 should obtain on at least $990$ of the $1000$ estimates $\widehat{\eta_{cv}}(q)$; in other words, at least $990$ of the $1000$ discrepancies $| \eta_{or} - \widehat{\eta_{cv}}(q)|$ should be bounded by
 \begin{equation}
 0 \leq | \widehat{\eta_{cv}}(q) - \eta_{or} | \leq \epsilon \eta_{or}.
 \label{theory+}
 \end{equation}
Using the relation $\eqref{rach}$ between $\epsilon$ and $r$, this bound becomes tighter as the number $r$ of CV measurements increases; however, at the same time, the oracle error $\eta_{or}$ increases with $r$ for fixed $m$ as fewer measurements $n = m - r$ are input to OMP.  An ideal number $r$ of CV measurements should not be too large or too small;  Figure 1 suggests that setting aside just enough measurements $r$ such that $\epsilon \leq.6$ is satisfied in $\eqref{rach}$ serves as a good heuristic to choose the number of cross validation measurements (in Figure 1, $\epsilon \leq .6$ is satisfied by taking only $r = 30$ measurements).  
 \\
 \\
We indicate the theoretical bound $\eqref{theory}$ with dark gray in Figure 1, which is compared to the interval in light gray of the $990$ values of $\eta_{cv}(q)$ that are closest to $\eta_{or}$ in actuality.   
 \\
 \\
This experiment is run for several values of $r$ within the interval $[5, 90]$; the results are plotted in Figure 1(a), with the particular range $r \in [5, 30]$ blown up in Figure 1(b). 
\\
\\
We have also carried out this experiment with a smaller noise variance; i.e. $x_b = x_0 + {\cal N}_b$ is subject to additive noise
\begin{equation}
{\cal N}_b \sim N(0, .02).
\end{equation}
The signal $x_b$ is again renormalized to satisfy $||x_b||_{l_2^N} = 1$; it now has an expected noise level of 
\begin{equation}
E(\sigma_d(x_b)) \approx .116.
\end{equation}
The results of this experiment are plotted in Figure 1(c).
\\
\\
\begin{figure}[htp]
{ \label{}
\includegraphics[width=6.5in]{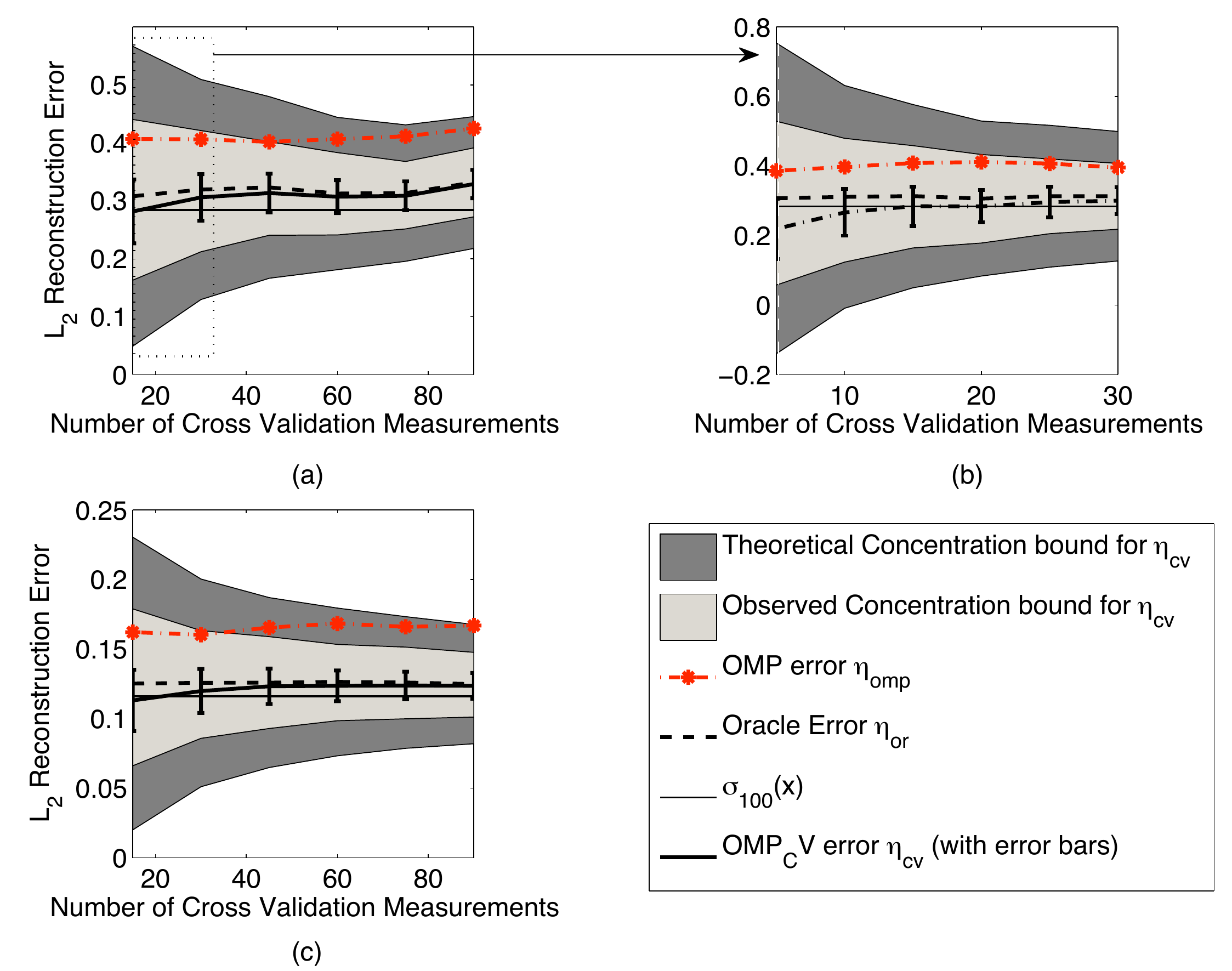}}

\caption{\small{Comparison of the reconstruction algorithms OMP and OMP-CV.  We fix the parameters $N = 3600, m = 800, k = 200$, and underlying sparsity $d = 100$, but vary the number $r$ of the total $m$ measurements reserved for cross validation from $5$ to $90$, using the remaining $n = m-r$ measurements for training.  The underlying signal has residual $\sigma_{100}(x) \approx .284$ in Figures 1(a) and 1(b), and $\sigma_{100}(x) \approx .116$ in Figure 1(c), as shown for reference by the thin horizontal line.  In both cases, the OMP-CV  error $\eta_{cv}$ (the solid black line with error bars; each point represents the average of 1000 trials) gives a better approximation to the residual error than does OMP (dot-dashed line) with very high probability, even when as few as $20$ of the total $800$ measurements are used for cross validation.  Even though $\eta_{cv}$ is guaranteed to provide a tighter bound for $\eta_{or}$ as the number $r$ of CV measurements increases, at the same time, the oracle error $\eta_{or}$ becomes a worse indicator of the residual $\sigma_{100}(x)$ because fewer measurements $n = m - r$ are input to OMP.  } }
\label{figure 1}
\end{figure}

\subsection{Experimental Results}
\begin{enumerate}

\item We remind the reader that the cross-validation estimates $\widehat{\eta_{cv}}$ are observable to the user, while the values of $\eta_{omp}$, $\eta_{or}$, along with the noise level $\sigma_d(x)$, are not available to the user.  Nevertheless, $\widehat{\eta_{cv}}$ can serve as a proxy for $\eta_{or}$ according to $\eqref{theory}$, and this is verified by the plots in Figure 1.  $\widehat{\eta_{cv}}$ can also provide an upper bound on $\sigma_d(x)$, as is detailed in Section 5.1.

\item The theoretical bound $\eqref{theory}$ is seen to be tight, when compared with the observed concentration bounds in Figure 1.  

\item With high probability, the estimate $\hat{x}_{cv}(15)$ using $r = 15$ out of the alloted $m = 800$ measurements will be a better estimate of $x$ than the OMP estimate: $|| \hat{x}_{cv}(15) - x ||_2 \leq || \hat{x}_{omp}(15) - x ||_2$.  With \emph{overwhelming probability}, the estimate $\hat{x}_{cv}(30)$ will result in error $|| \hat{x}_{cv}(30) - x ||_2 \leq || \hat{x}_{omp}(30) - x||_2$.  We note that the estimates $\hat{x}_{cv}(15)$ and $\hat{x}_{cv}(30)$ correspond to accuracy parameters $\epsilon(15) = .8405$ and $\epsilon(30) = .5943$ in \eqref{rach}, indicating that $\epsilon \leq.6$ is a good heuristic to determine when enough CV measurements have been reserved.

\item The OMP-CV estimate $\hat{x}_{cv}$ will have more pronounced improvement over the OMP estimate $\hat{x}_{omp}$ when there is larger discrepancy between the true sparsity $d$ of $x_0$ and the upper bound $k$ used by OMP (in Figure (1), $d = 100$ and $k = 200$).   In contrast, OMP-CV will not outperform OMP in approximation accuracy when $d$ is close to $k$; however, the multiplicative relation $\eqref{theory}$ guarantees that OMP-CV will not underperform OMP, either.   
\end{enumerate}

\section{Beyond Compressed Sensing}

The Compressed Sensing setup can be viewed within the more general class of {\it underdetermined linear inverse problems}, in which $x \in \mathbb{R}^N$ is to be reconstructed from a known $m \times N$ underdetermined matrix ${\cal A}$ and lower dimensional vector $y = {\cal A} x$ using a decoding algorithm $\Delta: \mathbb{R}^m \rightarrow \mathbb{R}^N$; in this broader context, ${\cal A}$ is given to the user, but not necessarily {\it specified by} the user as in compressed sensing.  In many cases, a prior assumption of sparsity is imposed on $x$, and an iterative decoding algorithm such as LASSO $\eqref{lasso}$ will be used to reconstruct $x$ from $y$ $\cite{17}$.   If it is possible to take on the order of $r = \log{p}$ additional measurements of $x$ by an $r \times N$ matrix $\Psi$ satisfying the conditions of Lemma $\eqref{psi}$, then all of the analysis presented in this paper applies to this more general setting.  In particular, the error $||x - \hat{x}_j ||_{l_2^N}$ at up to $j \leq p$ successive approximations $\hat{x}_j$ of the decoding algorithm $\Delta$ may be bounded from below and above using the quantities $|| \Psi (x - \hat{x}_j) ||_{\ell_2^r}$, and the final approximation $\hat{x}$ to $x$ can be chosen from among the entire sequence of estimates $\hat{x}_j$ as outlined in Theorem $\eqref{mainthm}$;  an earlier estimate $\hat{x}_j$ may approximate $x$ better than a final estimate $\hat{x}_p$ which contains the artifacts of parameter overfitting occurring at later stages of iteration. 

\section{Extensions and Open Problems}
We have presented an alternative approach to compressed sensing in which a certain number $r$ of the $m$ allowed measurements of a signal $x \in \mathbb{R}$$^N$ are reserved to track the error in decoding by the remaining $m - r$ measurements, allowing us to choose a best approximation to $x$ in the metric of $\ell^N_2$ out of a sequence of $p$ estimates $(\hat{x}_j)_{j=1}^p$, and estimate the error between $x$ and its best approximation by a $k$-sparse vector, again with respect to the metric of $\ell^N_2$.   We detailed how the number $r$ of such measurements should be chosen in terms of desired accuracy $\epsilon$ of estimation, confidence level $\xi$ in the prediction, and number $p$ of decoding iterations to be measured; in general, $r = O(\log(p))$ measurements suffice.   Several important issues remain unresolved; we mention only a few below.
\begin{enumerate}
\item The cross validation technique promoted in this paper corresponds specifically to the technique of \emph{holdout cross validation} in statistics, where a data set is partitioned into a single training and cross validation set (as a rule of thumb, the cross validation set is usually taken to be less than or equal to a third of the size of the training set; in the the current paper, we have shown that the Johnson Lindenstrauss lemma provides a theoretical justification of how many, or, more precisely, how few, cross validation measurements are needed in the context of compressed sensing).  Other forms of cross validation, such as \emph{repeated random subsampling} cross validation or \emph{K-fold} cross validation, remain to be analyzed in the context of compressed sensing.  The former technique corresponds to repeated application of holdout cross validation, with  $r$ cross validation measurements out of the total $m$ measurements chosen by random selection at each application.  The results are then averaged (or otherwise combined) to produce a single estimation.   The latter technique, $K$-fold cross validation, also corresponds to repeated application of holdout cross validation.  In this case, the $m$ measurements are partitioned into $K$ subsets of equal size $r$, and cross-validation is repeated exactly $K$ times with each of the $K$ subsets of measurements used exactly once as the validation set. The $K$ results are again combined to produce a single estimation.   Although Theorem $\eqref{mainthm}$ does not directly apply to these cross validation models, the experimental results of Section 6 suggest that, equiped with an $m \times N$ matrix satisfying the requirements of Lemma $\eqref{psi}$, the application of $K$ fold cross validation to subsets of the measurements of size $r << m - r$ just large enough that $\epsilon > 0$ in Theorem $\eqref{mainthm}$ for fixed accuracy $\xi$ and constant $C = 1$ can be combined to accurately approximate the underlying signal with near certainty.  

\item It is not clear that the analysis in Theorem $\eqref{mainthm}$ can be extended to the  noisy compressed sensing model,
\begin{equation}
y = \Phi x + {\cal N},
\label{noise}
\end{equation}
where ${\cal N} \sim N(0,\sigma^2)$ is a Gaussian random variable that accounts for both noise and quantization error on the measurements $\Phi x$.  Because measurement noise and quantization error are unavoidable in any real-world sensing device, any proposed compressed sensing technique should extend to the model $\eqref{noise}$.   Indeed, cross validation is studied in $\cite{7}$ in this context as a stopping criterion for decoding algorithms of homotopy/greedy type, in the case that $x$ is truly sparse and ${\cal N}$ is Gaussian noise.   The experimental results in $\cite{7}$ indicate that cross validation works well in this setting, but it remains to provide theoretical justification of these results.

\item We have only considered cross validation over the metric of $\ell_2$.  However, the error $||x - \hat{x}||_{\ell_2^N}$, or \emph{root mean squared error}, is just one of several metrics used in image processing for analyzing the quality of a reconstruction $\hat{x} \in \mathbb{R}^N$ to a (known) image $x \in \mathbb{R}^N$.  In fact, the $\ell_1$ reconstruction error $||x - \hat{x} ||_{\ell_1^N}$ has been argued to outperform the root mean squared error as an indicator of reconstruction quality $\cite{23}$.  Unfortunately, Theorem $\eqref{mainthm}$ cannot be extended to the metric of $\ell_1$, as there exists no $\ell_1$ analog of the Johnson Lindenstrauss Lemma \cite{27}.  However, it remains to understand the extent to which cross validation in compressed sensing can be applied over a broader class of image reconstruction metrics, perhaps using more refined techniques than those considered in this paper.
\item Many more compressed sensing matrices than just those satisfying the requirements of Lemma $\eqref{psi}$ are observed \emph{empirically} to satisfy the Johnson Lindenstrauss Lemma;  in particular, a properly normalized $r \times N$ matrix obtained by selecting $r$ rows at random from the $N \times N$ discrete Fourier matrix is observed to satisfy $\eqref{psi}$ for number of measurements $r \geq \epsilon^{-2}\log{\frac{1}{2\delta}}$, and the empirical concentration bounds for these matrices appear to be \emph{indistinguishable} from those of (properly normalized) random Gaussian and Bernoulli matrices of the same dimension.   This suggests that cross validation should be considered as a general technique that can be applied to a set of $m$ compressed sensing measurements, the theoretical justification of which is a very interesting problem that we hope to pursue in the future.   

\end{enumerate}

\section*{Acknowledgment}

The author would like to thank Ingrid Daubechies and Albert Cohen for their insights and encouragement, without which this paper could not have been written.

\end{document}